\documentclass[10pt]{amsart}

\usepackage{amsfonts,amsmath,latexsym,amssymb} 
\usepackage[T2A]{fontenc}
\usepackage[utf8]{inputenc}
\usepackage[english]{babel}

\newtheorem{theorem}{Theorem}
\newtheorem{lemma}{Lemma}
\newtheorem{remark}{Remark}
\newtheorem{definition}{Definition}
\newtheorem{corollary}{Corollary}
\DeclareMathOperator*{\esssup}{{\rm ess\,sup}}
\newcommand{\Rd}{{\mathbb R}^d}
\newcommand{\RR}{\mathbb R}

\newcommand{\polarK}{K^\circ}

\allowdisplaybreaks
\usepackage[left=3cm,right=2cm,
    top=2.5cm,bottom=2.5cm,bindingoffset=0cm]{geometry}

\begin{document}
\title{ On  Landau -- Kolmogorov type inequalities for charges and their applications}

\author[V.~F.~Babenko]{Vladyslav Babenko}
\address{Department of Mathematical Analysis and Optimization, Oles Honchar Dnipro National University, Dnipro, Ukraine}

\author[V.~V.~Babenko]{Vira Babenko}
\address{Department of Mathematics and Computer Science, Drake University, Des Moines, USA}

\author[O.~V.~Kovalenko]{Oleg Kovalenko}
\address{Department of Mathematical Analysis and Optimization, Oles Honchar Dnipro National University, Dnipro, Ukraine}

\author[N.~V.~Parfinovych]{Nataliia Parfinovych}
\address{Department of Mathematical Analysis and Optimization, Oles Honchar Dnipro National University, Dnipro, Ukraine}
\begin{abstract}
In this article we prove sharp Landau--Kolmogorov type inequalities on a class of charges defined on Lebesgue measurable subsets of a cone in $\RR^d$, $d\geq 1$, that are absolutely continuous with respect to the Lebesgue measure. In addition we  solve the Stechkin problem of approximation of the Radon--Nikodym derivative of such charges by bounded operators and two related problems. As an application, we also solve these extremal problems on classes of essentially bounded functions $f$ such that their distributional partial derivative  $\frac{\partial ^d f}{\partial x_1\ldots\partial x_d}$ belongs to the Sobolev space $W^{1,\infty}$.
\end{abstract}
\maketitle





\section{Introduction}
Inequalities that estimate norms of  intermediate derivatives of univariate or multivariate functions using the norms of the functions and their derivatives of higher order play an important role in many branches of Analysis and its applications. It appears that the richest applications are obtained from sharp inequalities of this kind, which attracts much interest to the inequalities with the smallest possible constants. For univariate functions, the results by Landau~\cite{Landau}, and Kolmogorov~\cite{Kolmogorov1939}, are among the brightest ones in this topic. Inequalities of this kind are often called the Landau-Kolmogorov type inequalities. A survey on the results for univariate and multivariate functions for the case of derivatives of integer order and further references can be found in~\cite{Arestov1996, BKKP, Konovalov1978,
BuslaevTikhomirov1979,Timoshin1995,
Timofeev1985,BKP1997,Babenko2000}. 

Let $X$ and $Y$ be linear spaces equipped with a seminorm $\|\cdot\|_X$ and a norm $\|\cdot\|_Y$ respectively. A linear operator $A\colon X\to Y$ is called bounded if 
\[
\| A\|=\| A\|_{X\to Y}:=\sup\limits_{\| x\|_X\le 1}\| Ax\|_Y<\infty.
\]
Otherwise, the operator $A$ is called unbounded.
By $\mathcal{L}(X,Y)$ we denote  the space of all linear bounded operators $S\colon X\to Y$. 

Let $A\colon X \to Y$ be an operator (not necessarily linear) with a domain $D_A \subset X$. Let also $\mathfrak{M}\subset D_A$ be some class of elements. The problem to find the modulus of continuity of the operator $A$ on the class $\mathfrak{M}$ i.e., the function 
\begin{equation}\label{modCont}
\Omega(\delta)=\Omega(A,\mathfrak{M};\delta):=\sup\{\|Ax\|_{Y}:x\in \mathfrak{M}, \ \|x\|_{X}\leq \delta\},\,\delta\geq 0,
\end{equation}
is an abstract version of the problem about the Landau -- Kolmogorov type inequalities.

For the operator  $A$ and an operator $S\in \mathcal{L}(X,Y)$ we set
$$
    U(A,S; \mathfrak{M}):=\sup\{\|Ax-Sx\|_Y\colon x\in \mathfrak{M}\} .
$$

The Stechkin problem of approximation of a generally speaking unbounded operator $A$ by linear bounded operators is stated as follows. For a given number $N> 0$ find the quantity
\begin{equation}\label{bestappr}
    E_N(A,\mathfrak{M}):=\inf\left\{ U(A,S;\mathfrak{M})\colon S\in{\mathcal L}(X,Y), \|S\|\le N\right\},
\end{equation}  
and an operator $S$ on which the infimum is attained, if such an operator exists.
The statement of this problem, first important 
 results, and
solutions to this problem for differential operators of small orders were presented in~\cite{Stechkin1967}. For a survey of further results on this problem see~\cite{ Arestov1996}. 

The Stechkin problem, in turn, is intimately connected to Landau-Kolmogorov type inequalities. The following theorem describes this connection, see~\cite{Stechkin1967}. 

\begin{theorem}\label{th::Stechkin1} For any $x\in \mathfrak{M}$ and arbitrary $S\in \mathcal{L}(X,Y)$ the following inequality holds
$$
 \| Ax\|\le  \| Ax-Sx\|+\| S\|\|x\|
\le U(A,S;\mathfrak{M})+\|S\|\cdot \|x\|,    
$$
    and, therefore,
$$
     \forall x\in \mathfrak{M},\;\forall N>0, \quad \|Ax\|\le E_N(A,\mathfrak{M})+N\|x\|.  
$$
Moreover
\[
\Omega(A,\mathfrak{M},\delta)\le \inf\limits_{N>0}\{E_N(A,\mathfrak{M})+N\delta\}.
\]
If in addition there exist $\overline{S}\in \mathcal{L}(X,Y)$ and  $\overline{x}\in Q$ such that
\begin{equation}\label{i::equality}
    \|A\overline{x}\|=\| A\overline{x}-\overline{S}\overline{x}\|+\| \overline{S}\|\|\overline{x}\|=U(A,\overline{S};\mathfrak{M})+\|\overline{S}\|\cdot \|\overline{x}\|,
\end{equation}
then 
$$
\Omega(\|\overline{x}\|_X)=\| A\overline{x}\|_Y\text{ and }E_{\left\|\overline{S}\right\|}(A, \mathfrak{M})=U(A,\overline{S};\mathfrak{M})=\| A\overline{x}\|-\left\|\overline{S}\right\|\|\overline{x}\|.
$$
Thus the element $\overline{x}$ is extremal for  problem~\eqref{modCont}, and the operator $\overline{S}$ is optimal for problem~\eqref{bestappr}.
\end{theorem}

\begin{remark}
In Stechkin's article~\cite{Stechkin1967} it is assumed that  $X$ and $Y$ are Banach spaces. However, as it is easy to see, completeness, and even presence of a norm in $X$ is not necessary. It is sufficient to have a seminorm in $X$. Completeness of $Y$ is also not necessary.
\end{remark}
Let $\mathfrak{M}\subset D(A)$ and $S\in\mathcal{L}(X,Y)$. For $\delta\geq 0$ the value
$$
U_{\delta}(A,S,\mathfrak{M}):=\sup\{\|Ax-S\eta\|_{Y}:x\in \mathfrak{M},\eta\in X,\|x-\eta\|_{X}\le\delta\},
$$
if called the deviation of the operator $S\in\mathcal{L}(X,Y)$ from the operator  $A$ on the class $\mathfrak{M}$, where elements are known with error $\delta$. The problem of optimal recovery of the operator $A$ by linear operators on such  
class $\mathfrak{M}$ consists of finding the quantity
$$
\mathcal{E}_{\delta}(A,\mathfrak{M},\mathcal{L}(X,Y)):=\inf\limits_{ S\in\mathcal{L}(X,Y)}U_{\delta}(A,S,\mathfrak{M})
$$
and an operator $S$, on which the infimum on the right-hand side is attained, if it exists.

The connection of this problem with the Landau -- Kolmogorov type inequalities and the Stechkin problem is given by the following theorem, see~\cite[Theorem 2.1]{Arestov1996}.
\begin{theorem}\label{th::Stechkin2}
    If $\mathfrak{M}$ is a convex centrally-symmetric set and $A$ is a homogeneous operator, then for all 
    $\delta>0$
    \[
    \Omega(A,\mathfrak{M};\delta)\le \mathcal{E}_{\delta}(A,\mathfrak{M},\mathcal{L}(X,Y))\le \inf\limits_{N>0}\{ E_N(A,\mathfrak{M})+N\delta\}.
    \]
    If in addition there exist an element $\overline{x}\in \mathfrak{M}$ and an operator $\overline{S}$ that satisfy~\eqref{i::equality}, then
    \[
   \mathcal{E}_{\|\overline{x}\|}(A,\mathfrak{M},\mathcal{L}(X,Y))=\Omega(A,\mathfrak{M};\|\overline{x}\|)=\| A\overline{x}\|.
    \]
\end{theorem}

We consider charges $\nu$ defined on Lebesgue measurable subsets of an open cone  $C\subset \RR^d$ that are absolutely continuous with respect to the Lebesgue measure $\mu$. For such charges we obtain sharp Landau -- Kolmogorov type inequalities that estimate the $L_\infty (C)$--norm  of the Radon -- Nikodym derivative $D_\mu\nu$ using the value of some seminorm of the charge $\nu$ and the $L_\infty (C)$--norm of the  gradient $\nabla D_\mu\nu$ of this derivative.  We also solve the problem of  approximation of the operator  $D_\mu$ by bounded operators on the class of charges $\nu$ such that the   $L_\infty (C)$--norm of the vector-function $\nabla D_\mu\nu$ does not exceed $1$, and the problem of optimal recovery of the operator $D_\mu$ on this class given the elements known with error. 

In the case $C=\RR_+^m\times \RR^{d-m}$, $0\le m\le d$, we obtain  inequalities  that estimate the $L_\infty(C)$--norm of a mixed derivative of the function $f\colon C\to \RR$ using the $L_\infty(C)$--norm of the function and the $L_\infty(C)$--norm of the gradient of the mixed derivative, and prove sharpness of the obtained inequalities in the cases $m=0,1$. We also solve the problem of approximation of the mixed derivative operator by bounded operators on the class that is defined by restrictions on the $L_\infty(C)$--norm of the gradient of the mixed derivative; and the problem of optimal recovery of the mixed derivative operator on this class with elements known with error.

The article is organized as follows. In Section~\ref{s::definitions} we give necessary definitions and state the main results of the paper. In Sections~\ref{s::L-KInequalities} and~\ref{s::mixedDerivative} we give the proofs for the stated results.

\section{Definitions and main results}\label{s::definitions}

For $x,y\in\Rd$, $d\geq 1$, denote by $(x,y)$ the dot product of $x$ and $y$. Let $K\subset \Rd$ be an open convex symmetric with respect to $\theta$ bounded set with $\theta\in{\rm int\,} K$.  For $x\in \Rd$ denote by 
$|x|_K$ the norm of $x$ generated by the set $K$ i.e.,
$|x|_K:=\inf\{\lambda>0\colon x\in \lambda K\}.$
Denote by $\polarK$ the polar set of the set $K$ i.e., 
$$K^\circ = \left\{y\in \Rd\colon \sup\limits_{x\in K} (x,y)\leq 1\right\}.$$
It is well known that the norm $|\cdot|_{\polarK}$ is the dual to $|\cdot|_K$ norm i.e., 
$
|z|_{\polarK} = \sup \{(x,z)\colon |x|_K\leq 1\}.
$ 
In particular for all $x,y\in\Rd$,
$
(x,y)\leq |x|_K\cdot |y|_{\polarK}.
$
Everywhere below $C\subset \RR^d$ is an open convex cone and $\mu$ is the Lebesgue measure in $\RR^d$.

By  $\mathfrak{N}(C)$ we denote the set of  charges $\nu$ (see e.g.,~\cite[Chapter 5]{Berezanski}) defined on Lebesgue measurable subsets of the open cone  $C\subset \RR^d$ that are absolutely continuous with respect to the Lebesgue measure $\mu$. By the Radon -- Nikodym theorem, there exists a finite integrable function $f\colon C\to \RR$ such that for each measurable set $Q\subset C$
\[
\nu(Q)=\int\limits_Qf(x)dx.
\]
This function $f$ is called the Radon--Nikodym derivative of $\nu$. We denote the Radom--Nikodym derivative of the measure $\nu$ by $D_\mu\nu$ i.e., $D_\mu\nu(x) = f(x)$ $\mu$--almost everywhere. 

The set $\mathfrak{N}(C)$ is a linear space with respect to the standard addition and multiplication by a real number. Consider a family of seminorms  $\{\|\cdot\|_{K,h}, h>0\}$ and a seminorm $\|\cdot\|_{K}$ in this space as follows:
\[
\|\nu\|_{K,h}=\|\nu(\cdot +hK\cap C)\|_{L_\infty(C)}\text { and } 
\|\nu\|_K:=\sup\limits_{h>0}\|\nu\|_{K,h}.
\]

Everywhere below, for a locally integrable function $f$ by $\nabla f$ we denote the gradient of the function $f$; the partial derivatives  are understood in the distributional sense.

\begin{definition}
    Let $Q\subset \Rd$ be an open set. By  $W^{1,\infty}(Q)$ we denote the Sobolev space of functions $f\colon Q\to \RR$ such that $f$ and all their (distributional) partial derivatives of the first order belong to $L_\infty(Q)$.
\end{definition}
\begin{remark}
    Due to equivalence of norms in finite dimensional spaces, for each $f\in W^{1,\infty}(Q)$ one has $|\nabla f|_{\polarK} \in L_\infty(Q)$.
\end{remark}
\begin{definition}
Denote by $\mathfrak{N}_K(C)$ ($\mathfrak{N}_{K, h}(C)$, $h>0$)  the set of all charges $\nu\in\mathfrak{N}(C)$ that have a finite seminorm  $\|\nu\|_K$ (resp. $\|\nu\|_{K,h}$), and by $\mathfrak{N}^{1,\infty}_K(C)$ ($\mathfrak{N}^{1,\infty}_{K,h}(C)$) the set of all charges $\nu \in \mathfrak{N}_K(C)$ (resp.  $\nu \in \mathfrak{N}_{K,h}(C)$)  such that $D_\mu\nu\in W^{1,\infty}(C)$.
 \end{definition} 

 We prove the following theorem, which gives sharp Landau -- Kolmogorov type inequalities in both additive and multiplicative forms for charges.

\begin{theorem}\label{th::LKInequalityForMeasures}
Let $h>0$ and $\nu\in \mathfrak{N}^{1,\infty}_{K,h}(C)$. Then $D_\mu\nu\in L_\infty(C)$ and the following additive Landau -- Kolmogorov inequality holds:
       \begin{equation}\label{LKadd}
   \left\| D_\mu\nu\right\|_{L_\infty(C)}
   \le
   \left\| D_\mu\nu-S_h\nu\right\|_{L_\infty(C)}+\| S_h\|\|\nu\|_{K,h}
   \\ \le 
    \frac {dh}{d+1}\left\|\left|\nabla D_\mu\nu\right|_{\polarK}\right\|_{L_\infty(C)}+\frac{\|\nu\|_{K,h}}{h^d\mu(K\cap C)},
\end{equation}
where 
\begin{equation}\label{extremalOperator}
S_h\nu(x)=\frac{\nu(x+hK\cap C)}{h^d\mu(K\cap C)}.
\end{equation}

    If $\nu\in \mathfrak{N}^{1,\infty}_K(C)$, then inequality~\eqref{LKadd} holds for any $h>0$, and the following multiplicative Landau -- Kolmogorov inequality holds: \begin{equation}\label{LKmult}
    \left\| D_\mu\nu\right\|_{L_\infty(C)}\le\left(\frac{d+1}{\mu(K\cap C)}\right)^{\frac 1{d+1}}\left\|\left|\nabla D_\mu\nu\right|_{\polarK}\right\|_{L_\infty(C)}^{\frac d{d+1}}\|\nu\|_K^{\frac 1{d+1}}.
\end{equation}
Both inequalities are sharp. They become equalities on measures 
\begin{equation}\label{extremalNu}
    Q\mapsto \int\limits_Q(h-|x|_K)_+(x)d\mu(x),
\end{equation}
where for $\alpha\in \RR$, $\alpha_+ = \max\{\alpha,0\}$; moreover, inequality~\eqref{LKmult} becomes equality on  measures~\eqref{extremalNu} with arbitrary $h>0$.
\end{theorem}

In the following theorem we find the modulus of continuity and solve the Stechkin problem for the operator $D_\mu$. We also solve the problem of optimal recovery of the operator $D_\mu$ given the elements known with error. 
\begin{theorem}\label{th::relatedProblemsForMeasures}
 Denote by $\mathfrak{M}$ the family of charges $\nu\in \mathfrak{N}^{1,\infty}_K(C)$ such that $\|\nabla D_\mu\nu\|_{L_\infty(C)}\le 1$. For all $\delta > 0$ and $N> 0$ the following equalities hold:
\begin{equation}\label{chargeModnepr}
\mathcal{E}_\delta(D_\mu,\mathfrak{M},\mathcal{L}) = \Omega(D_\mu,\mathfrak{M};\delta)=
    \left(\frac{(d+1)\delta}{\mu(K\cap C)}\right)^{\frac{1}{d+1}};
    \end{equation}
\begin{equation}\label{chargeAppr}
E_N(D_\mu,\mathfrak{M})=\frac d{d+1}\left(\frac{1}{N\mu(K\cap C)}\right)^{\frac{1}{d}}.
\end{equation}
\end{theorem}

 Set ${\bf I}=(1,\ldots,1)\in \RR^d$, and for a locally integrable function $f$
$$\partial_{\bf I}f(x)=\frac{\partial ^d f}{\partial x_1\ldots\partial x_d}(x),$$
where the derivatives are understood in the distributional sense. 

\begin{definition}
By $W^{\bf{I},\nabla}_{\infty,\infty}(C)$ we denote the space of functions  $f\in L_\infty(C)$ such that  $\partial_{\bf I}f\in W^{1,\infty}(C)$. 
\end{definition}

The following theorems holds.

\begin{theorem}\label{th::LKInequalitiesForDerivative}
Let $K=(-1,1)^d$ and therefore $|x|_K=\max\{|x_1|,\ldots,|x_d|\}$, and $|x|_{K^o}=\sum_{i=1}^d|x_i|$. Let also $C=\RR^d_{m,+}:=\RR^m_+\times\RR^{d-m}, 0\le m\le d$ and  $f\in W^{{\bf I},\nabla}_{\infty,\infty}(C)$. Then $\partial_{\bf I}f\in L_\infty(C)$ and the following inequalities hold.
\begin{enumerate} 
         \item The Landau -- Kolmogorov inequality in the additive form: 
      \begin{equation}\label{LKadd1}
     \forall h>0,\quad\| \partial_{\bf I}f\|_{L_\infty(C)}\le h\frac d{d+1}\left\|| \nabla \partial_{\bf I}f|_{\polarK}\right\|_{L_\infty(C)}+\frac {2^m}{h^d}\| f\|_{L_\infty(C)}.
\end{equation}
   \item The Landau -- Kolmogorov inequality in the multiplicative form:
   \begin{equation}\label{LKmult_1}
      \| \partial_{\bf I}f\|_{L_\infty(C)}\le \left(2^m(d+1)\| f\|_{L_\infty(C)}\right)^{\frac 1{d+1}}\left\|| \nabla \partial_{\bf I}f|_{\polarK}\right\|_{L_\infty(C)}^{\frac d{d+1}}.
\end{equation}
\end{enumerate}
For $m=0$ and $m = 1$ the inequalities are sharp.
\end{theorem}
\begin{remark}
It is not clear, whether the inequalities in Theorem~\ref{th::LKInequalitiesForDerivative} are sharp for $m=2,\ldots,d$ even in the case $d=2$.
\end{remark}
\begin{theorem}\label{th::relatedProblemsForDerivative}
Let $C$ and $K$ be as in Theorem~\ref{th::LKInequalitiesForDerivative}. Denote by $\mathfrak{M}$ the family of functions $f\in W^{\bf{I},\nabla}_{\infty,\infty}(C)$ such that $\| |\nabla\partial_{\bf I}f|_{ \polarK}\|_{L_\infty(C)}\le 1.$
For all $\delta> 0, N>0$ and $m\in \{0,1\}$ the following equalities hold.
$$
     E_\delta(\partial_{\bf I},\mathfrak{M})
     =
     \Omega\left(\partial_{\bf I},\mathfrak{M};\delta\right)
     =
    \left(2^m(d+1)\delta\right)^{\frac 1{d+1}};
$$
$$
E_N(\partial_{\bf I},\mathfrak{M})
=
\frac d{d+1}\left(\frac{2^m}{N}\right)^{\frac{1}{d}}.
$$
\end{theorem}

\section{Inequalities of Landau -- Kolmogorov type for charges and related problems}\label{s::L-KInequalities}
\subsection{Auxiliary results}
\begin{lemma}\label{l::|u|integral}
The following equality holds:
    \[
\int\limits_{hK\cap C}|u|_Kd\mu(u)=\frac{dh^{d+1}}{d+1}\mu(K\cap C).
\]
\end{lemma}
\begin{proof}
Using the layer cake representation, see e.g.,~\cite[Theorem 1.13]{lieb2001}, we obtain
\begin{multline*}
\int\limits_{hK\cap C}|u|_Kd\mu(u)
= 
\int_0^\infty \mu\{u\in hK\cap C\colon |u|_K > t\}dt 
\\= 
\int\limits_0^h(\mu(hK\cap C)-\mu(tK\cap C))dt
 =
\mu(K\cap C)\int\limits_0^h(h^d-t^d)dt
=
\frac{dh^{d+1}}{d+1}\mu(K\cap C).
\end{multline*}
\end{proof}

We need the following result, which is related to~\cite[Lemmas~4 and~5]{Babenko22}.
\begin{lemma}\label{l::nablaF_eNorm}
    For $h>0$ and $x\in C$ set $f_{e,h}(x) = (h - |x|_K)_+$. Then $f$ is continuous and $\left\|\left|\nabla f_{e,h}\right|_{\polarK}\right\|_{L_\infty(C)}=1.$   
\end{lemma}
\begin{proof}
Continuity of $f$ is obvious. We prove the remaining statement of the lemma in the case when $K$ is a polytope (i.e., a convex hull of a finite number of points) first. If $\gamma$ is a face of $K$, $\delta$ is the distance from $\theta$ to $\gamma$,  and $n=(n_1,\ldots,n_d)$ is the exterior unit normal of $\gamma$, then for arbitrary $y\in {\rm Int}\, ({\rm conv} \,(\{\theta\}\cup \gamma))$ one has $|y|_K=\frac{1}\delta (y,n)$. Hence
\[
\frac{\partial |\cdot|_K}{\partial x_j}(y)=\lim\limits_{\eta\to 0}\frac{| y+\eta e_j|_K-| y|_K}{\eta}=\lim\limits_{\eta\to 0}\frac{(y+\eta e_j,n)-(y,n)}{\eta\delta}
=\frac 1\delta (e_j,n)=\frac 1\delta n_j.
\]
Thus $\nabla |\cdot|_K(y)=\frac 1\delta n$. Since $|n|_{\polarK}=\delta$, see~\cite[Lemma~4]{Babenko22}, we obtain the required equality in the case, when $K$ is a polytope. The general case can be obtained by approximating the set $K$ by polytopes.
\end{proof}

\begin{lemma}\label{l::extremalCharge}
For each $h>0$ the charge $\nu_{e,h}$ defined in~\eqref{extremalNu} belongs to $\mathfrak{N}^{1,\infty}_K(C)$, and hence to  $\mathfrak{N}^{1,\infty}_{K,h}(C)$. Moreover,
$
\left\| D_\mu\nu_{e,h}\right\|_{L_\infty(C)}=h,
$ 
$
\left\|\left|\nabla D_\mu\nu_{e,h}\right|_{\polarK}\right\|_{L_\infty(C)}=1,
$
and 
$\|\nu_{e,h}\|_K=\|\nu_{e,h}\|_{K,h} = \frac{h^{d+1}}{d+1}\mu(K\cap C)$.
\end{lemma}
\begin{proof}
The first and the second equalities follow from the fact that $D_\mu\nu_{e,h} = f_{e,h}$ and Lemma~\ref{l::nablaF_eNorm}. Using Lemma~\ref{l::|u|integral}, for arbitrary $r>0$ we obtain
$$
\|\nu_{e,h}\|_{K,r} = \esssup_{y\in C} \int\limits_{y+rK\cap C}(h-|x|_K)_+d\mu(x)
\leq \int\limits_{hK\cap C}(h-|x|_K)d\mu(x)=\frac{h^{d+1}}{d+1}\mu(K\cap C).
$$
At the same time, the function $y\mapsto \nu_{e,h}(y+  hK\cap C)$ is continuous on $C$, hence $$\|\nu_{e,h}\|_{K,h}\geq \nu_{e,h}(hK\cap C) =\frac{h^{d+1}}{d+1}\mu(K\cap C),$$
and the last equality from the lemma follows.
\end{proof}

\begin{lemma}\label{l::operatorSh}
Let $\mathfrak{N}$ be $\mathfrak{N}_K(C)$ or $\mathfrak{N}_{K,h}(C)$, $h>0$. For the operator $S_h\colon \mathfrak{N}\to L_\infty (C)$ defined by~\eqref{extremalOperator}, one has 
\begin{equation}\label{normS}
   \|S_h\| =  \| S_h\|_{\mathfrak{N}\to L_\infty(C)}=\frac{1}{h^d\mu(K\cap C)}.
\end{equation} 
\end{lemma}
\begin{proof}
Since
\[
\| S_h\nu\|_{L_\infty(C)}
=
\frac{\|\nu(x+hK\cap C)\|_{L_\infty(C)}}{h^d\mu(K\cap C)}
\le
\frac{\|\nu\|_{K,h}}{h^d\mu(K\cap C)}
\leq 
\frac{\|\nu\|_{K}}{h^d\mu(K\cap C)},
\]
we obtain
$
   \| S_h\|_{\mathfrak{N}\to L_\infty(C)}\leq \frac{1}{h^d\mu(K\cap C)}.
$
On the other hand, due to Lemma~\ref{l::extremalCharge}, we obtain 
\[
\| S_h\|_{\mathfrak{N}\to L_\infty(C)}
\ge 
\frac{|S_h\nu_{e,h}(\theta)|}{\| \nu_{e,h}\|_{\mathfrak{N}}}
=
\frac{\nu_{e,h}(hK\cap C)}{h^d\mu(K\cap C)\| \nu_{e,h}\|_K}
=
\frac{1}{h^d\mu(K\cap C)},
\]
which implies the required.
\end{proof}

\begin{lemma}
    For each $\nu\in \mathfrak{N}^{1,\infty}_{K,h}(C)$, $h>0$, one has
    \begin{equation}\label{deviation}
\left\| D_\mu\nu-S_h\nu\right\|_{L_\infty(C)}\le \frac {dh}{d+1}\left\|\left|\nabla D_\mu\nu\right|_{\polarK}\right\|_{L_\infty(C)},
\end{equation}
where $S_h$ is defined by~\eqref{extremalOperator}. Inequality~\eqref{deviation} is sharp. It turns into equality for the charge $\nu_{e,h}$.
\end{lemma}
\begin{proof}
Using the fundamental theorem of calculus for distributions, see e.g.~\cite[Theorem 6.9]{lieb2001}, we obtain for almost all $x\in C$
\begin{gather*}
\left| D_\mu\nu(x)-S_h\nu(x)\right|=\frac{1}{h^d\mu(K\cap C)}\left|\;\int\limits_{x+hK\cap C}\left( D_\mu\nu(x)-D_\mu\nu(u)\right)d\mu(u)\right|
\\
=\frac 1{h^d\mu(K\cap C)}\left|\; \int\limits_{x+hK\cap C} \int\limits_0^1(u-x,\nabla D_\mu\nu (x+t(u-x)))dtd\mu(u)\right|
\\
\le \frac 1{h^d\mu(K\cap C)} \int\limits_{x+hK\cap C} \int\limits_0^1|u-x|_K\left|\nabla D_\mu\nu (x+t(u-x)))\right|_{\polarK}dtd\mu(u)
\\
\le \frac {\left\|\left|\nabla D_\mu\nu\right|_{\polarK}\right\|_{L_\infty(C)}}{h^d\mu(K\cap C)}\int\limits_{x+hK\cap C}|u-x|_Kd\mu(u)=\frac {\left\|\left|\nabla  D_\mu\nu\right|_{\polarK}\right\|_{L_\infty(C)}}{h^d\mu(K\cap C)}\int\limits_{hK\cap C}|u|_Kd\mu(u)
\\= 
\frac{dh}{d+1} \left\|\left|\nabla  D_\mu\nu\right|_{\polarK}\right\|_{L_\infty(C)},
\end{gather*}
 which implies~\eqref{deviation}. Since for the charges $\nu_{e,h}$ 
$$
\left\| D_\mu\nu_{e,h}-S_h\nu_{e,h}\right\|_{L_\infty(C)}\ge \left| D_\mu\nu_{e,h}(\theta)-S_h\nu_{e,h}(\theta)\right|
=h-\frac 1{h^d\mu(K\cap C)}\int\limits_{hK\cap C}(h-|u|_K)d\mu(u)=\frac {dh}{d+1},
$$
we obtain 
\[
\left\| D_\mu\nu_{e,h}-S_h\nu_{e,h}\right\|_{L_\infty(C)}\ge\frac {dh}{d+1}.
\]
Due to Lemma~\ref{l::extremalCharge}, $\left\|\nabla D_\mu\nu_{e,h}\right\|_{L_\infty(C)}=1$, hence inequality~\eqref{deviation} is sharp and becomes equality for $\nu =\nu_{e,h}$.
\end{proof}

\subsection{Proof of Theorem~\ref{th::LKInequalityForMeasures}}
Now we are ready to prove Theorem~\ref{th::LKInequalityForMeasures}.
\begin{proof}
Using the triangle inequality and inequalities~\eqref{deviation} and~\eqref{normS}, we obtain~\eqref{LKadd}.
Putting
\[
h=\left(\frac{\|\nu\|_K}{\mu(K\cap C)}\cdot \frac{d+1}{\left\|\left|\nabla D_\mu\nu\right|_{\polarK}\right\|_{L_\infty(C)}}\right)^{\frac{1}{d+1}}
\]
to the right-hand side of~\eqref{LKadd} (with  $\|\nu\|_{K}$ instead of $\|\nu\|_{K,h}$ on the right-hand side) we obtain inequality~\eqref{LKmult}.
Taking into account Lemma~\ref{l::extremalCharge}, we obtain sharpness of the inequalities.
\end{proof}
\begin{corollary}
For any function $f\in  L_1(C)\cap W^{1,\infty}(C)$ and any $h>0$ the following sharp Nagy type inequality holds:
\[
\left\| f\right\|_{L_\infty(C)}
    \le \frac {dh}{d+1}\left\|\left|\nabla f\right|_{\polarK}\right\|_{L_\infty(C)}+\frac{\| f\|_{L_1(C)}}{h^d\mu(K\cap C)}.
\]
\end{corollary}
\begin{proof}
    For any function $f\in  L_1(C)\cap W^{1,\infty}(C)$ the charge $\nu(Q) = \int_Q f(x)dx$ belongs to $\mathfrak{N}^{1,\infty}_K(C)$ and $\|\nu\|_K\leq \|f\|_{L_1(C)}$, thus the inequality follows from Theorem~\ref{th::LKInequalityForMeasures}. Moreover, for the extremal in Theorem~\ref{th::LKInequalityForMeasures} charge $\nu_{e,h}$, we have $\|\nu_{e,h}\|_K = \|D_\mu \nu_{e,h}\|_{L_1(C)}$, and hence the inequality becomes equality for the function $D_\mu \nu_{e,h}$. 
\end{proof}

\subsection{Proof of Theorem~\ref{th::relatedProblemsForMeasures}}
\begin{proof}
From Lemma~\ref{l::extremalCharge} it follows that for each $\delta> 0$ there exists $h>0$ such that $\|\nu_{e,h}\|_K = \delta$. Thus inequality~\eqref{LKmult} and the fact that it becomes equality on each charge $\nu_{e,h}$, $h>0$,  imply the second equality in~\eqref{chargeModnepr}.

Inequality~\eqref{LKadd} becomes equality for the charge $\nu_{e,h}$ for each $h>0$, thus for all $h>0$,

\begin{equation}\label{equality}
  \left\| D_\mu\nu_{e,h}\right\|_{L_\infty(C)}=\left\| D_\mu\nu_{e,h}-S_h\nu_{e,h}\right\|_{L_\infty(C)}+ \|S_h\|\|\nu_{e,h}\|_{K}.  
\end{equation}
Hence Theorem~\ref{th::Stechkin1} implies   
    \[
E_{\|S_h\|}\left(D_\mu,\mathfrak{M} \right)=\left\| D_\mu\nu_{e,h}-S_h\nu_{e,h}\right\|_{L_\infty(C)}=\frac {dh}{d+1},
    \]
which in view of Lemma~\ref{l::operatorSh} is equivalent to~\eqref{chargeAppr}. Finally, equality of the left-most and the right-most terms in~\eqref{chargeModnepr} follows from~\eqref{equality} and Theorem~\ref{th::Stechkin2}.
\end{proof}

\section{Inequalities that involve the gradient of mixed derivatives and their applications.}\label{s::mixedDerivative}
\subsection{Proof of Theorem~\ref{th::LKInequalitiesForDerivative}}
\begin{proof}
 For $C=\RR^d_{m,+}$ and $K=(-1,1)^d$, one has $|x|_K=\max\limits_{i=1,\ldots,d}|x_i|=| x|_{\infty}$, $|x|_{\polarK}=\sum\limits_{i=1}^d| x_i|=| x|_{1}$ and $\mu(hK\cap C)=2^{d-m}h^d.$
Let $\{e_i\}$ be the standard basis in $\RR^d$. For $i = 1,\ldots, d$, we set 
\[
\Delta^+_{i,h}f(x):=f(x+he_i)-f(x)\text{ and }
\Delta_{i,h}f(x):=f(x+he_i)-f(x-he_i).
\]
In virtue of the Fubini theorem, for almost all  $x\in C$ one has
\begin{equation}\label{mixdiff}
\int\limits_{x+hK\cap C} \partial_{\bf I}f(u)d\mu(u)=(\Delta^+_{1,h}\circ\ldots\circ\Delta^+_{m,h}\circ\Delta_{m+1,h}\circ\ldots\circ\Delta_{d,h})f(x).
\end{equation}

 For each $f\in W^{\bf{I},\nabla}_{\infty,\infty}(C)$ consider the following charge
$$
\nu_{\partial_{\bf I}f}(Q)=\int\limits_Q\partial_{\bf I}f(x)d\mu(x),\text{ so that } D_\mu\nu_{\partial_{\bf I}f}=\partial_{\bf I}f.
$$
Equality~\eqref{mixdiff} implies the following estimate for $\|\nu_{\partial_{\bf I}f}\|_{K,h}$ with arbitrary $h>0$.
\begin{equation}\label{est}
\|\nu_{\partial_{\bf I}f}\|_{K,h}=\sup\limits_{x\in C}\left|\;\int\limits_{x+{hK\cap C}} \partial_{\bf I}f(u)d\mu(u)\right|\le 2^d\| f\|_{L_\infty(C)},
\end{equation}
hence $\nu\in  \mathfrak{N}^{1,\infty}_K(C)$. Applying inequality~\eqref{LKadd} and taking into account estimate~\eqref{est}, we obtain
\begin{equation}\label{LKh1}
 \| \partial_{\bf I}f\|_{L_\infty(C)}
 \le
 \left\| \partial_{\bf I}f-S_h\nu_{\partial_{\bf I}f}\right\|_{L_\infty(C)}+ \| S_h\|\|\nu_
{\partial_{\bf I}f}\|_K 
 \\ \le 
 \frac {dh}{d+1}\left\|| \nabla \partial_{\bf I}f|_{\polarK}\right\|_{L_\infty(C)}
 +\frac {2^{m}}{h^d}\| f\|_{L_\infty(C)}.
 \end{equation}
This inequality, in particular, implies $\partial_{\bf I}f\in L_\infty(C)$. Define an operator $\overline{S}_h\colon L_\infty(C)\to L_\infty(C)$, setting
\[
\overline{S}_hf(x)=\frac 1{2^{d-m}h^d}\left(\Delta^+_{1,h}\circ\ldots\circ\Delta^+_{m,h}\circ\Delta_{m+1,h}\circ\ldots\circ\Delta_{d,h}\right)f(x), h > 0.
\]
It is clear that for each $f\in W^{\bf{I},\nabla}_{\infty,\infty}(C)$,
$
\overline{S}_hf(x)=(S_h\nu_{\partial_{\bf I}f})(x),
$
where the operator $S_h$ was defined by~\eqref{extremalOperator}.
It is also easy to see that $\left\|\overline{S}_h\right\|=\frac{2^m}{h^d}$. Inequality~\eqref{LKh1} can now be re-written as follows. 
\begin{equation}\label{LKh1modified}
    \| \partial_{\bf I}f\|_{L_\infty(C)}\le\left\| \partial_{\bf I}f-\overline{S}_hf\right\|_{L_\infty(C)}+\| \overline{S}_h\|\|f\|_{L_\infty(C)}\\
 \le \frac {dh}{d+1}\left\|| \nabla \partial_{\bf I}f|_{\polarK}\right\|_{L_\infty(C)}
 +\frac {2^{m}}{h^d}\| f\|_{L_\infty(C)},
\end{equation}
which implies inequality~\eqref{LKadd1}. Putting in this inequality
\[
 h=\left(\frac{2^m\| f\|_{L_\infty(C)}}{\left\|| \nabla \partial_{\bf I}f|_{\polarK}\right\|_{L_\infty(C)}}(d+1)\right)^{\frac 1{d+1}},
 \]
we obtain inequality~\eqref{LKmult_1}.

Next we show that the obtained inequalities are sharp in the cases $C=\RR^d$ and $C=\RR_+\times\RR^{d-1}$ i.e., for $m=0,1$. 

We start from the case $m=0$. Define
 \[
f(x)=\int\limits_0^{x_1}\ldots\int\limits_0^{x_d}f_{e,h}(u)du_d\ldots du_1,
 \]
 where $f_{e,h}$ is defined in Lemma~\ref{l::nablaF_eNorm}. Then  $\partial_{\bf I}f = f_{e,h}$
and hence
 $$
 \| \partial_{\bf I}f\|_{L_\infty(\RR^d)}=h,\text{ and } \left\|| \nabla \partial_{\bf I}f|_{\polarK}\right\|_{L_\infty(\RR^d)}=1.
 $$
 Moreover, applying Lemma~\ref{l::|u|integral} with $K = (-1,1)^d$ and $C = (0,\infty)^d$, we obtain
 \[
 \| f\|_{L_\infty(\RR^d)}=\int\limits_{[0,h]^d}(h-|x|_\infty)dx=\frac {h^{d+1}}{d+1}.
 \]
Direct computations show that for the function $f$ and $m = 0$, inequalities~\eqref{LKadd1},~\eqref{LKh1modified} and~\eqref{LKmult_1} become equalities.

Next we show that each of inequalities~\eqref{LKadd1},~\eqref{LKh1modified} and~\eqref{LKmult_1} are sharp in the case $m=1$ too.
There exists $0<a<h$ such that
\[
\int\limits_{\{ x\in hK\cap C\colon x_1<a\}}(h-|x|_K)dx=\int\limits_{\{ x\in hK\cap C\colon x_1>a\}}(h-|x|_K)dx=\frac 12\int\limits_{hK\cap C}(h-|x|_K)dx.
\]
The set $hK\cap C$ consists of $2^{d-1}$ equal cubes with edges' length equal to $h$ and such that one of the vertices of each of the cubes is $\theta=\theta_d$. The hyperplane  $x_1=a$ divides these cubes into pieces $c_1^-,\ldots, c_{2^{d-1}}^-$ that have $\theta$ among vertices, and pieces $c_1^+,\ldots ,c_{2^{d-1}}^+$ that have $(a,\theta_{d-1})$ among vertices. 
Then 
$$
\sum\limits_{i=1}^{2^{d-1}}\int\limits_{c_i^+}(h-|x|_K)dx
=
\int\limits_{\{ x\in hK\cap C\colon x_1>a\}}(h-|x|_K)dx
=
\int\limits_{\{ x\in hK\cap C\colon x_1<a\}}(h-|x|_K)dx
=
\sum\limits_{i=1}^{2^{d-1}}\int\limits_{c_i^-}(h-|x|_K)dx.
$$
From symmetry considerations it follows that the summands in the left-most part of the equality, and summands in the right-most part of the equality are equal. Hence for all $i,j=1,\ldots,2^{d-1}$
\[
\int\limits_{c_i^-}(h-|x|_K)dx=\int\limits_{c_j^+}(h-|x|_K)dx=\frac 1{2^d}\int\limits_{hK\cap C}(h-|x|_K)dx=\frac 1{2}\cdot \frac {h^{d+1}}{d+1}.
\]
We define an extremal in this case function as follows.
\[ g(x) =\int\limits_a^{x_1}\int\limits_0^{x_2}\ldots\int\limits_0^{x_d}f_{e,h}(u)du_d\ldots du_1,
 \]
 where $f_{e,h}$ is defined in Lemma~\ref{l::nablaF_eNorm}.
Then
 \[
 \| g\|_{L_\infty(C)}=\frac 1{2}\cdot\frac {h^{d+1}}{d+1},\,\| \partial_{\bf I}g\|_{L_\infty(C)}=h,\text{ and }\| |\nabla \partial_{\bf I}g|_{\polarK}\|_{L_\infty(C)}=1.
 \]

 Direct computations show that inequalities~\eqref{LKadd1},~\eqref{LKmult_1}, and~\eqref{LKh1modified} with $m=1$ become equalities on the function $g$.
\end{proof}

\subsection{Proof of Theorem~\ref{th::relatedProblemsForDerivative}}
Theorem~\ref{th::relatedProblemsForDerivative} can be received using Theorem~\ref{th::LKInequalitiesForDerivative} and the same arguments that were used to prove Theorem~\ref{th::relatedProblemsForMeasures} with the help of Theorem~\ref{th::LKInequalityForMeasures}.

\bibliographystyle{elsarticle-num}
\bibliography{mybibfile}

\end{document}